\documentstyle[12pt,epsf,amstex,graphics]{article}
\def\SBIMSMark#1#2#3{
 \font\SBF=cmss10 at 10 true pt
 \font\SBI=cmssi10 at 10 true pt
 \setbox0=\hbox{\SBF Stony Brook IMS Preprint \##1}
 \setbox2=\hbox to \wd0{\hfil \SBI #2}
 \setbox4=\hbox to \wd0{\hfil \SBI #3}
 \setbox6=\hbox to \wd0{\hss
             \vbox{\hsize=\wd0 \parskip=0pt \baselineskip=10 true pt
                   \copy0 \break%
                   \copy2 \break%
                   \copy4 \break}}
 \dimen0=\ht6   \advance\dimen0 by \vsize \advance\dimen0 by 8 true pt
                \advance\dimen0 by -\pagetotal
 \dimen2=\hsize \advance\dimen2 by .25 true in
  \openin2=publishd.tex
  \ifeof2\setbox0=\hbox to 0pt{}
  \else 
     \setbox0=\hbox to 3.1 true in{
                \vbox to \ht6{\hsize=3 true in \parskip=0pt  \noindent  
                {\SBI Published in modified form:}\hfil\break
                \input publishd.tex 
                \vfill}}
  \fi
  \closein2
  \ht0=0pt \dp0=0pt
 \ht6=0pt \dp6=0pt
 \setbox8=\vbox to \dimen0{\vfill \hbox to \dimen2{\copy0 \hss \copy6}}
 \ht8=0pt \dp8=0pt \wd8=0pt
 \copy8
 \message{*** Stony Brook IMS Preprint #1, #2 ***}
}

\def\temp{1.35}%
\let\tempp=\relax
\expandafter\ifx\csname psboxversion\endcsname\relax
\else
    \ifdim\temp cm>\psboxversion cm
    \else
      \let\temp=\psboxversion
      \let\tempp= 
    \fi
\fi
\tempp
\let\psboxversion=\temp
\catcode`\@=11
\def\psfortextures{
\def\PSspeci@l##1##2{%
\special{illustration ##1\space scaled ##2}%
}}%
\def\psfordvitops{
\def\PSspeci@l##1##2{%
\special{dvitops: import ##1\space \the\drawingwd \the\drawinght}%
}}%
\def\psfordvips{
\def\PSspeci@l##1##2{%
\d@my=0.1bp \d@mx=\drawingwd \divide\d@mx by\d@my
\includegraphics{##1\space}}}%
\def\psforoztex{
\def\PSspeci@l##1##2{%
\special{##1 \space
      ##2 1000 div dup scale
      \number-\psllx\space\space \number-\pslly\space\space translate
}}}%
\def\psfordvitps{
\def\dvitpsLiter@ldim##1{\dimen0=##1\relax
\special{dvitps: Literal "\number\dimen0\space"}}%
\def\PSspeci@l##1##2{%
\at(0bp;\drawinght){%
\special{dvitps: Include0 "psfig.psr"}
\dvitpsLiter@ldim{\drawingwd}%
\dvitpsLiter@ldim{\drawinght}%
\dvitpsLiter@ldim{\psllx bp}%
\dvitpsLiter@ldim{\pslly bp}%
\dvitpsLiter@ldim{\psurx bp}%
\dvitpsLiter@ldim{\psury bp}%
\special{dvitps: Literal "startTexFig"}%
\special{dvitps: Include1 "##1"}%
\special{dvitps: Literal "endTexFig"}%
}}}%
\def\psfordvialw{
\def\PSspeci@l##1##2{
\special{language "PostScript",
position = "bottom left",
literal "  \psllx\space \pslly\space translate
  ##2 1000 div dup scale
  -\psllx\space -\pslly\space translate",
include "##1"}
}}%
\def\psforptips{
\def\PSspeci@l##1##2{{
\d@mx=\psurx bp
\advance \d@mx by -\psllx bp
\divide \d@mx by 1000\multiply\d@mx by \xscale
\incm{\d@mx}
\let\tmpx\dimincm
\d@my=\psury bp
\advance \d@my by -\pslly bp
\divide \d@my by 1000\multiply\d@my by \xscale
\incm{\d@my}
\let\tmpy\dimincm
\d@mx=-\psllx bp
\divide \d@mx by 1000\multiply\d@mx by \xscale
\d@my=-\pslly bp
\divide \d@my by 1000\multiply\d@my by \xscale
\at(\d@mx;\d@my){\special{ps:##1 x=\tmpx cm, y=\tmpy cm}}
}}}%
\def\psonlyboxes{
\def\PSspeci@l##1##2{%
\at(0cm;0cm){\boxit{\vbox to\drawinght
  {\vss\hbox to\drawingwd{\at(0cm;0cm){\hbox{({\tt##1})}}\hss}}}}
}}%
\def\psloc@lerr#1{%
\let\savedPSspeci@l=\PSspeci@l%
\def\PSspeci@l##1##2{%
\at(0cm;0cm){\boxit{\vbox to\drawinght
  {\vss\hbox to\drawingwd{\at(0cm;0cm){\hbox{({\tt##1}) #1}}\hss}}}}
\let\PSspeci@l=\savedPSspeci@l
}}%
\newread\pst@mpin
\newdimen\drawinght\newdimen\drawingwd
\newdimen\psxoffset\newdimen\psyoffset
\newbox\drawingBox
\newcount\xscale \newcount\yscale \newdimen\pscm\pscm=1cm
\newdimen\d@mx \newdimen\d@my
\newdimen\pswdincr \newdimen\pshtincr
\let\ps@nnotation=\relax
{\catcode`\|=0 |catcode`|\=12 |catcode`|
|catcode`#=12 |catcode`*=14
|xdef|backslashother{\}*
|xdef|percentother{
|xdef|tildeother{~}*
|xdef|sharpother{#}*
}%
\def\R@moveMeaningHeader#1:->{}%
\def\uncatcode#1{%
\edef#1{\expandafter\R@moveMeaningHeader\meaning#1}}%
\def\execute#1{#1}
\def\psm@keother#1{\catcode`#112\relax}
\def\executeinspecs#1{%
\execute{\begingroup\let\do\psm@keother\dospecials\catcode`\^^M=9#1\endgroup}}%
\def\@mpty{}%
\def\matchexpin#1#2{
  \fi%
  \edef\tmpb{{#2}}%
  \expandafter\makem@tchtmp\tmpb%
  \edef\tmpa{#1}\edef\tmpb{#2}%
  \expandafter\expandafter\expandafter\m@tchtmp\expandafter\tmpa\tmpb\endm@tch%
  \if\match%
}%
\def\matchin#1#2{%
  \fi%
  \makem@tchtmp{#2}%
  \m@tchtmp#1#2\endm@tch%
  \if\match%
}%
\def\makem@tchtmp#1{\def\m@tchtmp##1#1##2\endm@tch{%
  \def\tmpa{##1}\def\tmpb{##2}\let\m@tchtmp=\relax%
  \ifx\tmpb\@mpty\def\match{YN}%
  \else\def\match{YY}\fi%
}}%
\def\incm#1{{\psxoffset=1cm\d@my=#1
 \d@mx=\d@my
  \divide\d@mx by \psxoffset
  \xdef\dimincm{\number\d@mx.}
  \advance\d@my by -\number\d@mx cm
  \multiply\d@my by 100
 \d@mx=\d@my
  \divide\d@mx by \psxoffset
  \edef\dimincm{\dimincm\number\d@mx}
  \advance\d@my by -\number\d@mx cm
  \multiply\d@my by 100
 \d@mx=\d@my
  \divide\d@mx by \psxoffset
  \xdef\dimincm{\dimincm\number\d@mx}
}}%
\newif\ifNotB@undingBox
\newhelp\PShelp{Proceed: you'll have a 5cm square blank box instead of
your graphics.}%
\def\s@tsize#1 #2 #3 #4\@ndsize{
  \def\psllx{#1}\def\pslly{#2}%
  \def\psurx{#3}\def\psury{#4}
  \ifx\psurx\@mpty\NotB@undingBoxtrue
  \else
    \drawinght=#4bp\advance\drawinght by-#2bp
    \drawingwd=#3bp\advance\drawingwd by-#1bp
  \fi
  }%
\def\sc@nBBline#1:#2\@ndBBline{\edef\p@rameter{#1}\edef\v@lue{#2}}%
\def\g@bblefirstblank#1#2:{\ifx#1 \else#1\fi#2}%
{\catcode`\%=12
\xdef\B@undingBox{
\def\ReadPSize#1{
 \readfilename#1\relax
 \let\PSfilename=\lastreadfilename
 \openin\pst@mpin=#1\relax
 \ifeof\pst@mpin \errhelp=\PShelp
   \errmessage{I haven't found your postscript file (\PSfilename)}%
   \psloc@lerr{was not found}%
   \s@tsize 0 0 142 142\@ndsize
   \closein\pst@mpin
 \else
   \if\matchexpin{\GlobalInputList}{, \lastreadfilename}%
   \else\xdef\GlobalInputList{\GlobalInputList, \lastreadfilename}%
     \immediate\write\psbj@inaux{\lastreadfilename,}%
   \fi%
   \loop
     \executeinspecs{\catcode`\ =10\global\read\pst@mpin to\n@xtline}%
     \ifeof\pst@mpin
       \errhelp=\PShelp
       \errmessage{(\PSfilename) is not an Encapsulated PostScript File:
           I could not find any \B@undingBox: line.}%
       \edef\v@lue{0 0 142 142:}%
       \psloc@lerr{is not an EPSFile}%
       \NotB@undingBoxfalse
     \else
       \expandafter\sc@nBBline\n@xtline:\@ndBBline
       \ifx\p@rameter\B@undingBox\NotB@undingBoxfalse
         \edef\t@mp{%
           \expandafter\g@bblefirstblank\v@lue\space\space\space}%
         \expandafter\s@tsize\t@mp\@ndsize
       \else\NotB@undingBoxtrue
       \fi
     \fi
   \ifNotB@undingBox\repeat
   \closein\pst@mpin
 \fi
\message{#1}%
}%
\def\psboxto(#1;#2)#3{\vbox{%
   \ReadPSize{#3}%
   \advance\pswdincr by \drawingwd
   \advance\pshtincr by \drawinght
   \divide\pswdincr by 1000
   \divide\pshtincr by 1000
   \d@mx=#1
   \ifdim\d@mx=0pt\xscale=1000
         \else \xscale=\d@mx \divide \xscale by \pswdincr\fi
   \d@my=#2
   \ifdim\d@my=0pt\yscale=1000
         \else \yscale=\d@my \divide \yscale by \pshtincr\fi
   \ifnum\yscale=1000
         \else\ifnum\xscale=1000\xscale=\yscale
                    \else\ifnum\yscale<\xscale\xscale=\yscale\fi
              \fi
   \fi
   \divide\drawingwd by1000 \multiply\drawingwd by\xscale
   \divide\drawinght by1000 \multiply\drawinght by\xscale
   \divide\psxoffset by1000 \multiply\psxoffset by\xscale
   \divide\psyoffset by1000 \multiply\psyoffset by\xscale
   \global\divide\pscm by 1000
   \global\multiply\pscm by\xscale
   \multiply\pswdincr by\xscale \multiply\pshtincr by\xscale
   \ifdim\d@mx=0pt\d@mx=\pswdincr\fi
   \ifdim\d@my=0pt\d@my=\pshtincr\fi
   \message{scaled \the\xscale}%
 \hbox to\d@mx{\hss\vbox to\d@my{\vss
   \global\setbox\drawingBox=\hbox to 0pt{\kern\psxoffset\vbox to 0pt{%
      \kern-\psyoffset
      \PSspeci@l{\PSfilename}{\the\xscale}%
      \vss}\hss\ps@nnotation}%
   \global\wd\drawingBox=\the\pswdincr
   \global\ht\drawingBox=\the\pshtincr
   \global\drawingwd=\pswdincr
   \global\drawinght=\pshtincr
   \baselineskip=0pt
   \copy\drawingBox
 \vss}\hss}%
  \global\psxoffset=0pt
  \global\psyoffset=0pt
  \global\pswdincr=0pt
  \global\pshtincr=0pt 
  \global\pscm=1cm 
}}%
\def\psboxscaled#1#2{\vbox{%
  \ReadPSize{#2}%
  \xscale=#1
  \message{scaled \the\xscale}%
  \divide\pswdincr by 1000 \multiply\pswdincr by \xscale
  \divide\pshtincr by 1000 \multiply\pshtincr by \xscale
  \divide\psxoffset by1000 \multiply\psxoffset by\xscale
  \divide\psyoffset by1000 \multiply\psyoffset by\xscale
  \divide\drawingwd by1000 \multiply\drawingwd by\xscale
  \divide\drawinght by1000 \multiply\drawinght by\xscale
  \global\divide\pscm by 1000
  \global\multiply\pscm by\xscale
  \global\setbox\drawingBox=\hbox to 0pt{\kern\psxoffset\vbox to 0pt{%
     \kern-\psyoffset
     \PSspeci@l{\PSfilename}{\the\xscale}%
     \vss}\hss\ps@nnotation}%
  \advance\pswdincr by \drawingwd
  \advance\pshtincr by \drawinght
  \global\wd\drawingBox=\the\pswdincr
  \global\ht\drawingBox=\the\pshtincr
  \global\drawingwd=\pswdincr
  \global\drawinght=\pshtincr
  \baselineskip=0pt
  \copy\drawingBox
  \global\psxoffset=0pt
  \global\psyoffset=0pt
  \global\pswdincr=0pt
  \global\pshtincr=0pt 
  \global\pscm=1cm
}}%
\def\psbox#1{\psboxscaled{1000}{#1}}%
\newif\ifn@teof\n@teoftrue
\newif\ifc@ntrolline
\newif\ifmatch
\newread\j@insplitin
\newwrite\j@insplitout
\newwrite\psbj@inaux
\immediate\openout\psbj@inaux=psbjoin.aux
\immediate\write\psbj@inaux{\string\joinfiles}%
\immediate\write\psbj@inaux{\jobname,}%
\def\toother#1{\ifcat\relax#1\else\expandafter%
  \toother@ux\meaning#1\endtoother@ux\fi}%
\def\toother@ux#1 #2#3\endtoother@ux{\def\tmp{#3}%
  \ifx\tmp\@mpty\def\tmp{#2}\let\next=\relax%
  \else\def\next{\toother@ux#2#3\endtoother@ux}\fi%
\next}%
\let\readfilenamehook=\relax
\def\re@d{\expandafter\re@daux}
\def\re@daux{\futurelet\nextchar\stopre@dtest}%
\def\re@dnext{\xdef\lastreadfilename{\lastreadfilename\nextchar}%
  \afterassignment\re@d\let\nextchar}%
\def\stopre@d{\egroup\readfilenamehook}%
\def\stopre@dtest{%
  \ifcat\nextchar\relax\let\nextread\stopre@d
  \else
    \ifcat\nextchar\space\def\nextread{%
      \afterassignment\stopre@d\chardef\nextchar=`}%
    \else\let\nextread=\re@dnext
      \toother\nextchar
      \edef\nextchar{\tmp}%
    \fi
  \fi\nextread}%
\def\readfilename{\bgroup%
  \let\\=\backslashother \let\%=\percentother \let\~=\tildeother
  \let\#=\sharpother \xdef\lastreadfilename{}%
  \re@d}%
\xdef\GlobalInputList{\jobname}%
\def\psnewinput{%
  \def\readfilenamehook{
    \if\matchexpin{\GlobalInputList}{, \lastreadfilename}%
    \else\xdef\GlobalInputList{\GlobalInputList, \lastreadfilename}%
      \immediate\write\psbj@inaux{\lastreadfilename,}%
    \fi%
    \let\readfilenamehook=\relax%
    \ps@ldinput\lastreadfilename\relax%
  }\readfilename%
}%
\def\nowarnopenout{%
 \def\warnopenout##1##2{%
   \readfilename##2\relax
   \message{\lastreadfilename}%
   \immediate\openout##1=\lastreadfilename\relax}}%
\def\warnopenout#1#2{%
 \readfilename#2\relax
 \def\t@mp{TrashMe,psbjoin.aux,psbjoint.tex,}\uncatcode\t@mp
 \if\matchexpin{\t@mp}{\lastreadfilename,}%
 \else
   \immediate\openin\pst@mpin=\lastreadfilename\relax
   \ifeof\pst@mpin
     \else
     \edef\tmp{{If the content of this file is precious to you, this
is your last chance to abort (ie press x or e) and rename it before
retexing (\jobname). If you're sure there's no file
(\lastreadfilename) in the directory of (\jobname), then go on: I'm
simply worried because you have another (\lastreadfilename) in some
directory I'm looking in for inputs...}}%
     \errhelp=\tmp
     \errmessage{I may be about to replace your file named \lastreadfilename}%
   \fi
   \immediate\closein\pst@mpin
 \fi
 \message{\lastreadfilename}%
 \immediate\openout#1=\lastreadfilename\relax}%
{\catcode`\%=12\catcode`\*=14
\gdef\splitfile#1{*
 \readfilename#1\relax
 \immediate\openin\j@insplitin=\lastreadfilename\relax
 \ifeof\j@insplitin
   \message{! I couldn't find and split \lastreadfilename!}*
 \else
   \immediate\openout\j@insplitout=TrashMe
   \message{< Splitting \lastreadfilename\space into}*
   \loop
     \ifeof\j@insplitin
       \immediate\closein\j@insplitin\n@teoffalse
     \else
       \n@teoftrue
       \executeinspecs{\global\read\j@insplitin to\spl@tinline\expandafter
         \ch@ckbeginnewfile\spl@tinline
       \ifc@ntrolline
       \else
         \toks0=\expandafter{\spl@tinline}*
         \immediate\write\j@insplitout{\the\toks0}*
       \fi
     \fi
   \ifn@teof\repeat
   \immediate\closeout\j@insplitout
 \fi\message{>}*
}*
\gdef\ch@ckbeginnewfile#1
 \def\t@mp{#1}*
 \ifx\@mpty\t@mp
   \def\t@mp{#3}*
   \ifx\@mpty\t@mp
     \global\c@ntrollinefalse
   \else
     \immediate\closeout\j@insplitout
     \warnopenout\j@insplitout{#2}*
     \global\c@ntrollinetrue
   \fi
 \else
   \global\c@ntrollinefalse
 \fi}*
\gdef\joinfiles#1\into#2{*
 \message{< Joining following files into}*
 \warnopenout\j@insplitout{#2}*
 \message{:}*
 {*
 \edef\w@##1{\immediate\write\j@insplitout{##1}}*
\w@{
\w@{
\w@{
\w@{
\w@{
\w@{
\w@{
\w@{
\w@{
\w@{
\w@{\string\input\space psbox.tex}*
\w@{\string\splitfile{\string\jobname}}*
\w@{\string\let\string\autojoin=\string\relax}*
}*
 \expandafter\tre@tfilelist#1, \endtre@t
 \immediate\closeout\j@insplitout
 \message{>}*
}*
\gdef\tre@tfilelist#1, #2\endtre@t{*
 \readfilename#1\relax
 \ifx\@mpty\lastreadfilename
 \else
   \immediate\openin\j@insplitin=\lastreadfilename\relax
   \ifeof\j@insplitin
     \errmessage{I couldn't find file \lastreadfilename}*
   \else
     \message{\lastreadfilename}*
     \immediate\write\j@insplitout{
     \executeinspecs{\global\read\j@insplitin to\oldj@ininline}*
     \loop
       \ifeof\j@insplitin\immediate\closein\j@insplitin\n@teoffalse
       \else\n@teoftrue
         \executeinspecs{\global\read\j@insplitin to\j@ininline}*
         \toks0=\expandafter{\oldj@ininline}*
         \let\oldj@ininline=\j@ininline
         \immediate\write\j@insplitout{\the\toks0}*
       \fi
     \ifn@teof
     \repeat
   \immediate\closein\j@insplitin
   \fi
   \tre@tfilelist#2, \endtre@t
 \fi}*
}%
\def\autojoin{%
 \immediate\write\psbj@inaux{\string\into{psbjoint.tex}}%
 \immediate\closeout\psbj@inaux
 \expandafter\joinfiles\GlobalInputList\into{psbjoint.tex}%
}%
\def\centinsert#1{\midinsert\line{\hss#1\hss}\endinsert}%
\def\psannotate#1#2{\vbox{%
  \def\ps@nnotation{#2\global\let\ps@nnotation=\relax}#1}}%
\def\pscaption#1#2{\vbox{%
   \setbox\drawingBox=#1
   \copy\drawingBox
   \vskip\baselineskip
   \vbox{\hsize=\wd\drawingBox\setbox0=\hbox{#2}%
     \ifdim\wd0>\hsize
       \noindent\unhbox0\tolerance=5000
    \else\centerline{\box0}%
    \fi
}}}%
\def\at(#1;#2)#3{\setbox0=\hbox{#3}\ht0=0pt\dp0=0pt
  \rlap{\kern#1\vbox to0pt{\kern-#2\box0\vss}}}%
\newdimen\gridht \newdimen\gridwd
\def\gridfill(#1;#2){%
  \setbox0=\hbox to 1\pscm
  {\vrule height1\pscm width.4pt\leaders\hrule\hfill}%
  \gridht=#1
  \divide\gridht by \ht0
  \multiply\gridht by \ht0
  \gridwd=#2
  \divide\gridwd by \wd0
  \multiply\gridwd by \wd0
  \advance \gridwd by \wd0
  \vbox to \gridht{\leaders\hbox to\gridwd{\leaders\box0\hfill}\vfill}}%
\def\fillinggrid{\at(0cm;0cm){\vbox{%
  \gridfill(\drawinght;\drawingwd)}}}%
\def\textleftof#1:{%
  \setbox1=#1
  \setbox0=\vbox\bgroup
    \advance\hsize by -\wd1 \advance\hsize by -2em}%
\def\textrightof#1:{%
  \setbox0=#1
  \setbox1=\vbox\bgroup
    \advance\hsize by -\wd0 \advance\hsize by -2em}%
\def\endtext{%
  \egroup
  \hbox to \hsize{\valign{\vfil##\vfil\cr%
\box0\cr%
\noalign{\hss}\box1\cr}}}%
\def\frameit#1#2#3{\hbox{\vrule width#1\vbox{%
  \hrule height#1\vskip#2\hbox{\hskip#2\vbox{#3}\hskip#2}%
        \vskip#2\hrule height#1}\vrule width#1}}%
\def\boxit#1{\frameit{0.4pt}{0pt}{#1}}%
\catcode`\@=12 
\psfordvips   

\newtheorem{PROP}{\bf Proposition}
\newtheorem{DEF}{\bf Definition}
\newtheorem{THEO}{\bf Theorem}
\newtheorem{LEMME}{\bf Lemma}
\newtheorem{CORO}{\bf Corollary}
\newtheorem{CONJ}{\bf Conjecture}

\def\Cm{{{\Bbb C}}}
\def\Bm{{{\Bbb B}}}
\def\Tm{{{\Bbb T}}}
\def\Zm{{{\Bbb Z}}}
\def\Dm{{{\Bbb D}}}
\def\Rm{{{\Bbb R}}}
\def\Nm{{{\Bbb N}}}
\def\Qm{{{\Bbb Q}}}
\def\M{{{\Bbb M}}}
\def\P{{{\Bbb P}}}
\def\Hm{{{\Bbb H}}}
\def\l{\ell}
\def\dem{\noindent{\sc Proof. }}
\def\findem{\hfill{\hbox{%
  \hskip 1pt%
  \vrule width 7pt height 6pt depth 1.5pt%
  \hskip 1pt}}}
\def\finlem{\hfill{$\square$}}

\pagestyle{myheadings}
\markboth{\centerline{\sc Xavier Buff}}{\centerline{Geometry of the 
Feigenbaum map.}} 
\setlength\textheight{9in} 
\setlength\textwidth{6.5in} 
\setlength\oddsidemargin{0pt} 
\setlength\evensidemargin{0pt} 
\setlength\topmargin{0pt} 
\addtolength\topmargin{-\headheight}   
\addtolength\topmargin{-\headsep}

\begin{document}
\SBIMSMark{1997/18}{November 1997}{}
\thispagestyle{empty}
\vskip 2cm
\centerline{\Large Geometry of the Feigenbaum map.}
\centerline{by}
\centerline{\sc \large Xavier Buff}
\centerline{\bf Cornell University,}
\centerline{100 White Hall,}
\centerline{Ithaca, NY, 14853}
\vskip .5cm
\centerline{\sc \Large Abstract.}

\vskip .5cm
We show that the Feigenbaum-Cvitanovi\'c equation can be interpreted 
as a linearizing equation, and the domain of 
analyticity of the Feigenbaum fixed point of renormalization as a 
basin of attraction. There is a natural decomposition of this basin 
which enables to recover a result of local connectivity by Jiang and 
Hu for the Feigenbaum Julia set.
\vskip 1cm
We have included in this article several computer drawn pictures. We would 
like to thank especially Henri Epstein who gave us the picture of 
the domain of analyticity of the Feigenbaum map, 
Louis Granboulan who helped us drawing some of the pictures using 
the program Pari, and Dan S{\o}rensen who designed a wonderful program to 
draw Mandelbrot and Julia sets.

I am grateful to Adrien Douady for having taught me so much about dynamical 
systems and renormalization. Most of the 
tools to deal with the Feigenbaum equation were explained to me by Henri 
Epstein, and most of my ideas came during discussions with him in IHES, with
Marguerite Flexor in Orsay and with John H. Hubbard in Cornell. I am greatly 
indebted to those four people.

\section*{Introduction.}
The notion of renormalization for dynamical systems was introduced by 
Feigenbaum and Cvitanovi\'c. Landford and Sullivan have then proved the 
existence of fixed points of renormalization. Those fixed points satisfy 
a functional equation which has been studied by Eckmann, Epstein, Wittwer 
and others. We will make use of this functional equation, known as the 
Feigenbaum-Cvitanovi\'c equation, to give a new approach to a result 
by Jiang and Hu regarding to the 
local connectivity of the Julia set of the Feigenbaum polynomial.

We assume that the reader is familiar with the notion of
renormalization and the notion of polynomial-like mappings.
One can for example read \cite{dh}, \cite{mcm1} or 
\cite{s}. The central object in our study is the Feigenbaum polynomial.
It is the most famous example of polynomial
which is infinitely
renormalizable (i.e. $k$-renormalizable for infinitely many
$k$). It is the unique real quadratic polynomial $z^2+c_{Feig}$
which is  
$2^k$-renormalizable for all $k\geq 1$.

\begin{figure}[htb]
\centerline{
\psboxscaled{500}{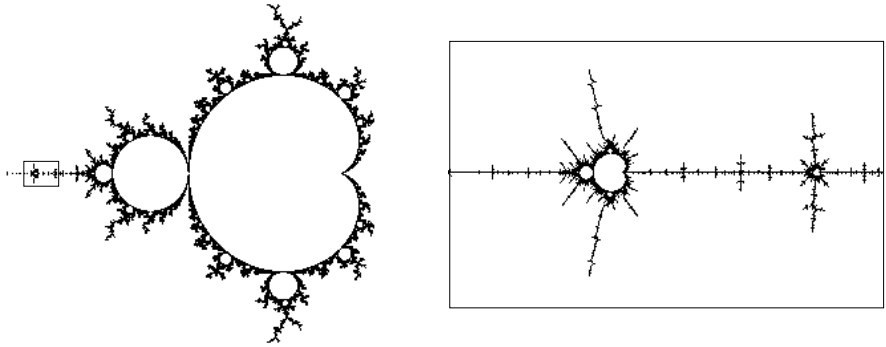}}
\caption{The Mandelbrot set and a small copy of it centered at a periodic 
point of period 3.}
\label{mandelbrot}
\end{figure}

One can define the Feigenbaum polynomial as the unique real polynomial
which is a fixed point of tuning by $-1$.
Tuning is the inverse of
renormalization. Given a parameter $c\in M$, such that $0$ is a periodic point 
of period $p$, Douady and Hubbard \cite{dh}
have constructed a tuning map, $x\rightarrow c*x$, which is a
homeomorphism of $M$ into itself, sending $0$ to $c$, and such that if
$x\neq 1/4$, then $f_{c*x}$ is $p$-renormalizable, and the
corresponding renormalization is in the same inner class as $f_x$ 
(where $f_x$ is the complex polynomial $f_x(z)=z^2+x$). 
This is the way they show there are small copies of the Mandelbrot set
inside itself (cf figure \ref{mandelbrot}).

The Feigenbaum value, $c_{Feig}=-1.401155...$, is in the intersection of all
the copies of $M$ obtained by tuning by $-1$. This intersection is not
known to consist of a single point, but its intersection with the real
axis is exactly the point $c_{Feig}$.

By construction, the Feigenbaum polynomial, $P_{Feig}$, is
$2$-renormalizable. 
There are several mappings $g:~U\rightarrow V$ such that
$g=f^2|U$ is a polynomial-like map with connected
Julia set.  All those maps have the same Julia set, and
are equal on this Julia set. We will say they define the same germ
$[g]$ of polynomial-like map (cf \cite{mcm1}).

We can define a renormalization operator ${\cal R}_2$ in the
following way. 
\begin{DEF}
Assume $[f]$ is a germ of a 2-renormalizable polynomial-like map. 
There exist open sets $U$ and $V$ such that the map $g:~U\to V$ 
defined by $g=f^2|U$ is a polynomial-like map with connected 
Julia set. The renormalization operator ${\cal R}_2$ is defined by
$${\cal R}_2([f])=[\alpha^{-1}\circ g \circ\alpha],$$
with $\alpha=g(0)=f^2(0)$, and $\alpha(z)=\alpha z$.
\end{DEF}
We have normalized the germ so that the critical value is 1.
 
\begin{DEF}
Two polynomial-like maps $f$ and $g$ representing germs are said to be
 hybrid equivalent if there 
is a quasi-conformal conjugacy between them with $\overline\partial h=0$ 
almost everywhere on the filled Julia set $K(f)$.
\end{DEF}
The equivalence classes are called inner classes.
One of the main theorems in \cite{dh} states that in the inner class any 
quadratic-like map with connected Julia set, there
is a unique polynomial $P(z)=z^2+c$, $c\in\Rm$.

Hence, if $[f]$ is a germ of Feigenbaum polynomial-like
map, then ${\cal R}_2([f])$ is hybrid equivalent to a unique polynomial. 
This polynomial is $2^k$-renormalizable for all $k\geq
1$. So it is the Feigenbaum 
polynomial. Hence ${\cal R}_2([f])$ is still a germ of Feigenbaum map, 
and we can iterate this
process, defining in such a way a sequence of germs: ${\cal
R}_2^n([P_{Feig}])$, $n\in\Nm$. The
following result has been proved and can be found in \cite{s} or 
\cite{mcm1}.
\begin{THEO}\label{fixpoint}
The sequence of germs ${\cal R}_2^n([P_{Feig}])$, $n\in\Nm$, converges
to a point $[\phi]$. This point is a fixed point of renormalization:
$${\cal R}_2([\phi])=[\phi],$$
and is in the inner class of the Feigenbaum polynomial. It is the 
unique fixed point of ${\cal R}_2$.
\end{THEO}

The properties we have mentioned still hold if the degree of the
critical point is any even integer. 
In the following $\l$ is an even integer
and $f$ is an analytic function which coincides in a neighborhood of
0 with the real fixed point of ${\cal R}_2$ with critical point of
degree $\l$.

In section 1, we state some results about this fixed point of 
renormalization. It is just a germ of polynomial-like map. However, 
we will show that there exists a natural representative of this germ. This 
construction has already been made by Epstein in \cite{e}, but we think it 
is worth making again here, to help the reader get accustomed with the 
tools.

In section 2, we give a description of the domain of analyticity of the 
fixed point of renormalization.  McMullen \cite{mcm1} proved that it has a 
maximal analytic extension $\hat f$ to a 
dense, simply connected open set $\widehat W\subset\Cm$. 
In lemmas 1 and 2, we prove that it is the basin of attraction of 
the map $\hat f(\lambda z)$, where $\lambda=-f^2(0)$. In particular, this 
enables us to prove that $\widehat W\subset \widehat W/\lambda$
(theorem~\ref{theoinclusion}),  
and to give a dynamical interpretation of the intersection $\displaystyle 
\bigcap_{n\in\Nm}\lambda^n\widehat W$ (theorem \ref{theointersection}). 

The boundary of $\widehat W$ is a closed set with empty interior. We would 
like to try and give a description of it. We prove that the boundary of 
$\widehat W$ contains accessible points 
with at least two accesses (proposition \ref{remark1}). This proves, in
particular, that  
this boundary does not have the structure of a Cantor Bouquet. Moreover, 
we think that both its Hausdorff dimension and its Lebesgue measure 
ought to be studied. The result of proposition \ref{remark2} is a step in 
this direction.

In section 3, we show that the domain $\widehat W$ is naturally paved by
puzzles that are homothetic to each other and cut the Julia set 
$K(f)$ in a connected way. This enables us to give a new proof of the local 
connectivity of the Julia set $K(f)$ at the critical point.

\section{Cvitanovi\'c-Feigenbaum equation.}

First of all, we would like to recall that $f$ is a solution of the following
system of equations (cf \cite{b2} or \cite{e}). 
\begin{DEF}
The Cvitanovi\'c-Feigenbaum equation:
$$\left\{\begin{array}{lll}
f(z) & = & -{1/\lambda}f\circ f(\lambda z),\\
f(0) & = & 1,\\
f(z) & = & F(z^\l),\mbox{ with }F^{-1}\mbox{ univalent in }
{\Bbb C}\setminus(]-\infty,-1/\lambda]\cup[1/\lambda^2,+\infty[).
\end{array}
\right.
$$
\end{DEF}

We will first state some results by Henri Epstein in the following two
propositions. Figure \ref{f} illustrates proposition \ref{epstein1}.
\begin{figure}[htb]
\centerline{
\input{f1.pstex_t}}
\caption{The graph of $f$ on $\Rm^+$.}
\label{f}
\end{figure}

\begin{figure}[htb]
\centerline{
\input{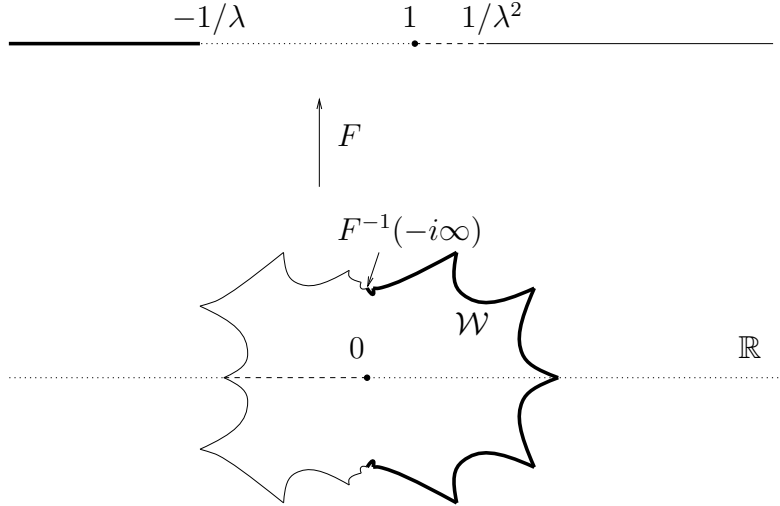}}
\caption{Maximal univalent extension of $F$.} 
\label{F}
\end{figure}

\begin{PROP} (cf \cite{e})\label{epstein1}
Let $f$ be a solution of the equation, and let $x_0$ be 
the first positive preimage of 0 by $f$. 
Then,
\begin{itemize}
\item[$\bullet$]{$f(1)=-\lambda$,}
\item[$\bullet$]{$f(\lambda x_0)=x_0$, and}
\item[$\bullet$]{the first critical point in $\Rm^+$ is $x_0/\lambda$,
with $f(x_0/\lambda)=-1/\lambda$.} 
\end{itemize}
\end{PROP}
This graph enables us to deduce the relative positions of some 
points on the real axis.
\begin{PROP} (cf \cite{e})\label{epstein2}
Univalent extension of $F$:
\begin{itemize}
\item[$\bullet$]{it is possible to extend $F^{-1}$ continuously to the 
boundary $\Rm$ of $\Hm_+$,
and even analytically except at points $(-1/\lambda)^n$, $n\geq
1$, which are branching points of type $z^{1/\l}$,}
\item[$\bullet$]{the values of $F^{-1}$ are never real except in
$[-1/\lambda,1/\lambda^2]$,}
\item[$\bullet$]{the extension of $F^{-1}$ to the closure of $\Hm_+$
is injective, and}
\item[$\bullet$]{when $z$ tends to infinity in $\Hm_+$,
$F^{-1}(z)$ tends to a point in $\Hm_-$
which will be denoted by $F^{-1}(i\infty)$.}
\end{itemize}
\end{PROP}
By symmetry, similar statements hold in $\Hm_-$.
Hence, $\cal W$ is a bounded domain of $\Cm$.
Those results are summarized in figure \ref{F}.

\begin{figure}
\vskip 12pt
\centerline{
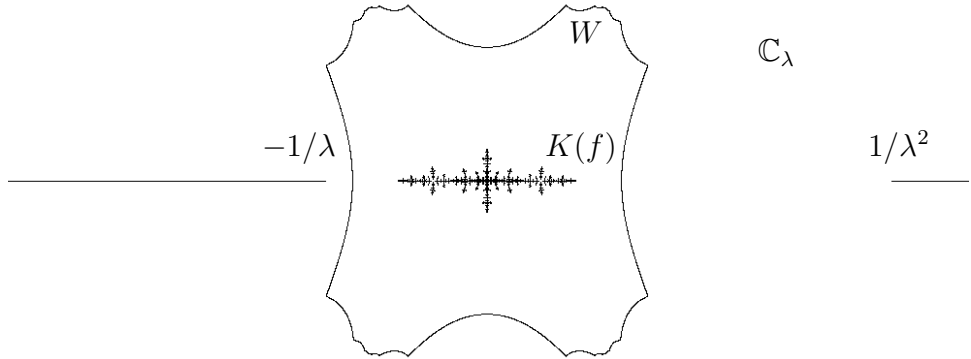}
\caption{A representative of the Feigenbaum fixed point of
renormalization, and the Julia set $K(f)$.}
\label{feigmap}
\end{figure}

In the following, we will use the notations:
\begin{eqnarray*}
{\Bbb C}_\lambda & = & {\Bbb C}\setminus(]-\infty,-{1\over\lambda}]
\cup[{1\over\lambda^2},+\infty[),\\
{\cal W} & = & F^{-1}(\Cm_\lambda), \\
W & = & \{z\in\Cm~|~z^\l\in {\cal W}\},\\
W_+ & = & W\cap\{z\in\Cm~|~0<\mbox{Arg}(z)<\pi/\l\},\mbox{ and}\\
W_- & = & W\cap\{z\in\Cm~|~-\pi/\l<\mbox{Arg}(z)<0\}.\\
\end{eqnarray*}
As $W_+\cap\Rm=[0,x_0/\lambda[$, we can deduce that
$W\cap \Rm=]-x_0/\lambda,x_0/\lambda[$ (here we use the 
fact that $\l$ is an even degree).
\begin{CORO}(cf figure \ref{feigmap})
The map $f~:~W\to \Cm_\lambda$ is a polynomial-like map of degree $\l$ 
representing the renormalization fixed point.
\end{CORO}
\dem
The graph of $f$ (cf figure \ref{f}) enables us to conclude that
$\overline{W}\subset \Cm_\lambda$, because $x_0<1$.
Hence $f~:~W\to \Cm_\lambda$ is a polynomial-like
map. Besides, $f$ has a unique critical point in 0 of degree $\l$. Then 
by uniqueness of the fixed point, the corollary follows.\findem

Hence, this polynomial-like map is quasi-conformally conjugate to the 
Feigenbaum polynomial. Thus, 
to know the geometry of the Julia set of the Feigenbaum
polynomial, it is enough to know the geometry of the Julia set
$K(f)$ of this polynomial-like map (cf figure \ref{feigmap}).

\section{Geometry of the domain of analyticity of $f$.}

Our first goal is to describe the geometry of the domain of analyticity 
of $f$. In \cite{mcm1}, McMullen proved that this domain exists, and is a 
dense open subset of $\Cm$. We will show that it is contained in one of 
its homothetics, and that it can be seen as the basin of attraction of a 
map related to $f$. This will enable us to introduce a chess board 
and some puzzles which we 
would like to use to prove the local connectivity of the Julia set $K(f)$.

\subsection{The domain of analyticity.}
\begin{DEF}
Let $f$ and $g$ be two holomorphic functions defined on open connected
domains of $\Cm$: $U_f$ and $U_g$. We say $g$ is an analytic extension
of $f$ if $f=g$ on some nonempty open set. 
Moreover, if all such analytic extensions are restrictions of a 
single map
$$\hat f:~\widehat W\rightarrow \Cm,$$
we will say that $\hat f$ is the unique maximal analytic extension 
of $f$
\end{DEF}

We will use a result by
McMullen concerning the domain of analyticity of the Feigenbaum map $f$.
\begin{PROP}(cf \cite{mcm1}).\label{mcmullen1}
Let $f$ be a solution of the Cvitanovi\'c-Fibonacci equation.
There exists a unique maximal extension of $f$,
$$\hat f:~\widehat W\rightarrow \Cm.$$
where $\widehat W$ is a simply connected open set in $\Cm$.
\end{PROP}
\dem 
First of all, recall that $f:~W\to \Cm_\lambda$ is a polynomial-like
map with non-escaping critical point. Hence, the set
$W_n=f^{-2^n}(\Cm_\lambda)$ is a simply connected open set and the map 
$f^{2^n}:~W_n\to \Cm_\lambda$ is a proper map. But the
Cvitanovi\'c-Feigenbaum equation enables us to say that
$$1/\lambda^n f^{2^n} (\lambda^n z):~W_n/\lambda^n\to \Cm_\lambda/\lambda^n$$
is an extension of $f$. Moreover, as $\Cm_\lambda/\lambda^n\subset
\Cm_\lambda/\lambda^{n+1}$, 
$$W_n/\lambda^n\subset W_{n+1}/\lambda^{n+1}.$$
We can now define 
$$\widehat W=\bigcup_{n\in\Nm} W_n/\lambda^n.$$
The map $f$ has a unique maximal analytic extension on $\widehat W$ which
coincides with $1/\lambda^n f^{2^n} (\lambda^n z)$, if $\lambda^n z\in
W_n$. \findem

\begin{figure}
\begin{center}
\begin{picture}(435,291)
\put(-1,-52){\includegraphics{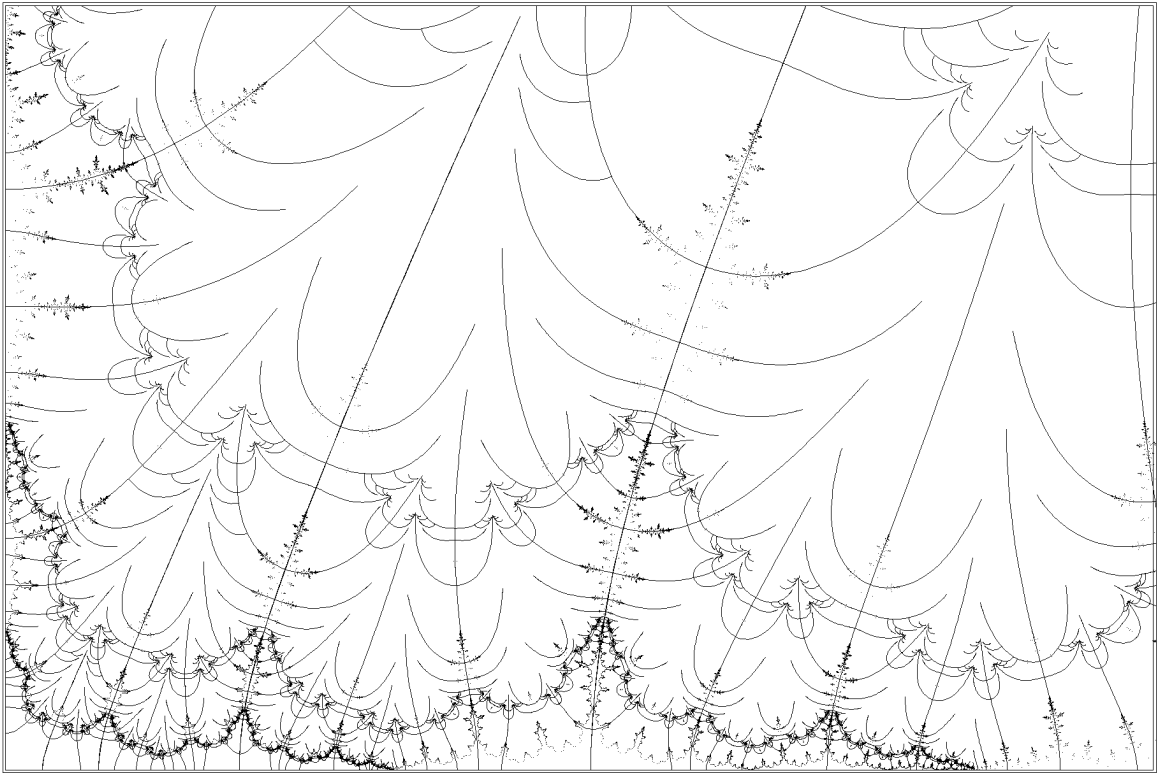}}
\end{picture}
\end{center}
\caption{The skeleton of the domain $\widehat W$, and some points in 
its boundary. The reader should imagine the domain by putting some 
flesh on the skeleton.}
\label{domaine}
\end{figure}

Figure~\ref{domaine} was given to us by Henri Epstein. This 
picture should enable the reader to imagine the 
intersection of the domain $\widehat W$ with $\{z=x+iy\in \Cm~|~x\geq 0~{\rm 
and}~ y\geq 0\}$. The lines are preimages of the real axis under $\hat f$ and 
form the skeleton of $\widehat W$. The dots are points in the boundary of 
$\widehat W$. The reader should imagine the domain $\widehat W$ by 
putting some flesh between the skeleton and the dots. But one should 
remember that the domain is simply connected, and dense in the 
plane. Hence, each point in the boundary of $\widehat W$ is connected 
to infinity.

\subsection{A linearizing equation.}
Now, we will show that the set $\widehat W$ can be seen as a basin of
attraction. 
\begin{LEMME}
The map $f_\lambda(z)=f(\lambda z)$ 
has an attracting fixed point $x_0$ of multiplier $-\lambda$. Moreover
$f_\lambda$ exchanges the two sets $W_+$ and $W_-$.
\end{LEMME}
\dem The trick in the proof is that the Cvitanovi\'c-Feigenbaum 
equation can be translated on the following commutative diagram: 
$$\begin{array}{rcl}
(\widehat W,x_0) & \stackrel{f(\lambda z)}{\longrightarrow} & (\widehat 
W,x_0)  \\ 
{\scriptstyle f}\downarrow &  & \downarrow {\scriptstyle f} \\
(\Cm,0) & \underset{-\lambda z}{\longrightarrow} & (\Cm,0). 
\end{array}
$$
It is a linearizing equation which proves the first part of the proposition. 
On the other hand,
$$\Hm_-=f(W_+)=-1/\lambda f(f(\lambda W_+)).$$
Hence $f(f_\lambda(W_+))=\Hm_+$. The set $f_\lambda(W_+)$ is a preimage 
of $\Hm_+$ under $f$. As
$x_0$ is in the closure of $W_+$ and $f_\lambda(x_0)=x_0
\in\overline{W_-}$, we can immediately deduce that $f_\lambda(W_+)=W_-$.
We can use the same arguments to show that $f_\lambda(W_-)=W_+$.

\findem

It is now possible to show the following result.
\begin{LEMME}\label{2}
The set $\widehat W$ is the basin of attraction of the fixed point $x_0$ of
$\hat f_\lambda=\hat f(\lambda z)$.
\end{LEMME}
\dem The commutative diagram we have written tells us that the linearizer
of the map $f_\lambda$ is $f$. But we know the domain of analyticity
of the linearizer is the immediate basin of the attracting
fixed point $x_0$. 
\finlem

We will prove that $\widehat W$ is dense in $\Cm$. Hence, there cannot be 
any other component in the basin of $\hat f_\lambda$.

\subsection{An inclusion of sets with their homothetics.}
To study the geometry of $\widehat W$, McMullen proved the following
result.
\begin{LEMME}\label{3}
The sets $K(f)/\lambda^n$, $(n\in\Nm)$, are all contained in $\widehat
W$.
\end{LEMME}
\dem The domain $\widehat W$ is the union of the sets
$W_n/\lambda^n$. But $K(f)\subset W_n=f^{-2^n}(\Cm_\lambda)$. Hence, for all
$n\in\Nm$,  
$$K(f)/\lambda^n\subset W_n/\lambda^n\subset \widehat W.$$
\finlem

\begin{THEO}\label{theoinclusion}
We have the following inclusion of sets.
$$\lambda K(f)\subset K(f)\subset
\bigcup_{n\in\Nm}\frac{K(f)}{\lambda^n} \subset \bigcap_{n\in\Nm}
\lambda^n \widehat W \subset \widehat W\subset \frac{\widehat
W}{\lambda}.$$ 
\end{THEO}
\dem The set $\lambda K(f)$ is the Julia set $K(f^2)$ of 
the renormalization 
$$f^2~:~\lambda W\to \lambda \Cm_\lambda.$$
 This Julia 
set is contained in $K(f)$. The inclusion  
$$K(f)\subset \bigcup_{n\in\Nm}K(f)/\lambda^n$$ 
follows immediately.
The inclusion $\widehat W\subset \widehat W/\lambda$ comes from lemma
\ref{2}. The domain $\widehat W$ is the basin of attraction of $\hat
f_\lambda$. Hence, it is contained in the domain of analyticity of
$\hat f_\lambda$ which is $\widehat W/\lambda$. The inclusion 
$$\bigcap_{n\in\Nm} \lambda^n \widehat W\subset \widehat W$$
follows immediately from this one. The remaining inclusion 
$$\bigcup_{n\in\Nm}\frac{K(f)}{\lambda^n} \subset \bigcap_{n\in\Nm}
\lambda^n \widehat W$$ comes from lemma \ref{3}. \findem

McMullen proved that
the union $\bigcup_{n\in\Nm}K(f)/\lambda^n$ is dense in $\Cm$. Hence
$\widehat W$ is a dense open subset of $\Cm$. As $K(f)$ is a compact
set with empty interior, the inclusions
are strict because of Baire Theorem: a 
countable union of closed sets with empty interior cannot be equal to
a countable intersection of dense open sets.

Finally, we would like to mention that the intersection of the
homothetics $\lambda^n \widehat W$ can be defined dynamically.
\begin{THEO}\label{theointersection}
The intersection of all the sets $\lambda^n \widehat W$ is the set of 
points the orbit of which under $\hat f$ stays in $\widehat W$:
$$\bigcap_{n\in\Nm} \lambda^n \widehat W = \widehat K=\{z\in\widehat
W~|~(\forall n\in\Nm)~\hat f^n(z)\in\widehat W\}.$$
\end{THEO}
\dem
To show this equality, note that 
$$\widehat K=\bigcap_{n\in\Nm} \hat f^{-n}(\widehat W).$$
We have noticed that $\hat f(\lambda \widehat W)=\hat
f_\lambda(\widehat W)=\widehat W$. 
Moreover, in a neighborhood of the origin, 
$$\hat f(z)=-\frac{1}{\lambda}\hat f\circ\hat f(\lambda z).$$
This equality has to be true whenever both sides are simultaneously
defined. The left side is defined on $\widehat W$. When $z\in\widehat
W$, $\lambda z\in\lambda\widehat W$. Hence, $f(\lambda z)\in \widehat
W$, and the right side of the equality is defined. The
Cvitanovi\'c-Feigenbaum equation holds on the whole domain $\widehat
W$. As $\hat f^{-1}(\widehat W)=\lambda \widehat W$, it follows that $\hat
f^{-2^n}(\widehat W)=\lambda^n \widehat W$, and the theorem holds.
\findem 

\subsection{Two remarks on the boundary of $\widehat W$.}

The domain $\widehat W$ is dense in $\Cm$, hence its boundary has empty 
interior. However, it seems to have a very complex structure.
Studying its topology, its Hausdorff dimension, 
or its Lebesgue measure could 
be interesting. But we do not have many tools at the moment. 
However, we will show that this 
boundary does not have the structure of a Cantor Bouquet, by proving that 
some points have at least two accesses in $\widehat W$.
Moreover, we will prove that some points in the 
boundary of $\widehat W$ cannot be deep points of $\partial \widehat W$ in 
the sense of McMullen (cf \cite{mcm1}).

First of all, recall that Proposition \ref{epstein2} tells that 
the following decomposition of the boundary of $\partial W_+$
is well defined.

\begin{DEF}
Let us denote by $\partial W_+$ the boundary of $W_+$. 
This boundary is the union of 
\begin{itemize}
\item{an arc $\gamma_1$ which is mapped by $\hat f$ to $\Rm_+$,}
\item{an arc $\gamma_2$ which is mapped by $\hat f$ to $\Rm-$, and}
\item{a point $x_1$ such that $(x_1)^\l=F^{-1}(-i\infty)$.}
\end{itemize}
\end{DEF}


\begin{PROP}\label{remark1}
The point $x_1/\lambda$ belongs to the boundary of $\widehat W$, and the arcs 
$\gamma_1/\lambda$ 
and $\gamma_2/\lambda$ which both land at $x_1/\lambda$
belong to the domain of analyticity $\widehat W$. Moreover, they do not 
belong to the same access to $x_1/\lambda$ in $\widehat W$.
\end{PROP}

\dem
First of all, note that a slight improvement of Lemma 1 enables us to 
show that $(x_1,\overline x_1)$ is a repelling cycle of period 2 for the map 
$f_\lambda$. In fact, this is the way Epstein proves the last statement of 
Proposition \ref{epstein2} (cf \cite{e}).

Hence, $x_1$ does not belong to the basin $\widehat W$. Thus, it is in 
the boundary $\partial \widehat W$. Moreover, according to theorem 
\ref{theoinclusion}, 
$\lambda \widehat W\subset \widehat W$, thus $x_1$ belongs to the boundary 
of $\lambda \widehat W$, 
and the statement concerning $x_1/\lambda$ is proved.

The Feigenbaum map $f$ maps $\gamma_1$ and $\gamma_2$ to intervals of $\Rm$. 
Hence, $\hat f_\lambda$ maps $\gamma_1/\lambda$ and $\gamma_2/\lambda$ 
inside $\widehat W$. Now, recall that $\widehat W$ is the basin of 
attraction of $\hat f_\lambda$. Hence, $\gamma_1/\lambda$ and 
$\gamma_2/\lambda$ are in the basin of attraction of $\hat f_\lambda$, 
which concludes the first statement relative to $\gamma_1$ and $\gamma_2$.

To conclude the proof, note that the union 
$\gamma_1/\lambda\cup\gamma_2/\lambda\cup\{x_1/\lambda\}$ 
bounds the set $W_+/\lambda$ which contains the point 
$x_1\in\Cm\setminus\widehat W$ in its interior.
\findem

We will finally prove the following proposition.
\begin{PROP}\label{remark2}
The point $x_1\in\Cm\setminus \widehat W$ cannot be a 
deep point of $\Cm\setminus\widehat W$.
\end{PROP}
This means that there exists a constant $K<1$ such that for all radii
$r\in\Rm_+$, there exists a ball of radius $Kr$ inside the intersection 
$\widehat W\cap B(x_1,r)$.

\dem
Recall that all the critical values of $f_\lambda$ are on the real axis. 
Hence one can choose a branch of $f_\lambda^{-2}$ which maps $\Hm_+$ into 
itself, having $x_1$ as attracting fixed point.

Finally, Lemma 2 proves that this inverse branch maps the boundary of 
$\widehat W$ into itself. As its interior is empty, a classical result 
enables us to conclude that $x_1$ cannot be a deep point of $\Cm\setminus 
\widehat W$.
\findem

The same result holds for all the points which 
eventually land under iteration of $\hat f_\lambda$ on $x_1$. Obviously, 
those points are in the boundary of the basin $\widehat W$ and they cannot 
be deep points of $\partial \widehat W$.

\section{Local connectivity of $K(f)$.}

We will now give a new proof of the result by Jiang and Hu \cite{jh}
relative to the local
connectivity of the Julia set of the Feigenbaum polynomial, by using 
the fixed point of renormalization.

We have proved that $\widehat W$ is the basin of attraction of 
$\hat f_\lambda$. Moreover $\hat f_\lambda$ exchanges the two sets $W_+$ and 
$W_-$. Hence it is possible to define a chessboard for $\hat f_\lambda$. 
\begin{DEF}The squares of the chessboard are defined in the following way:
\begin{itemize}
\item[$\bullet$]{The first two chess squares are the sets $W_+$ and
$W_-$.}
\item[$\bullet$]{ The other chess squares are the preimages of
those two by the map $\hat f_\lambda$.}
\end{itemize}
\end{DEF}
Some of those chess squares appear at the bottom left of figure 
\ref{domaine}.

It is now natural to take into account the following puzzles.

\begin{figure}[htb]
\centerline{
\psboxscaled{500}{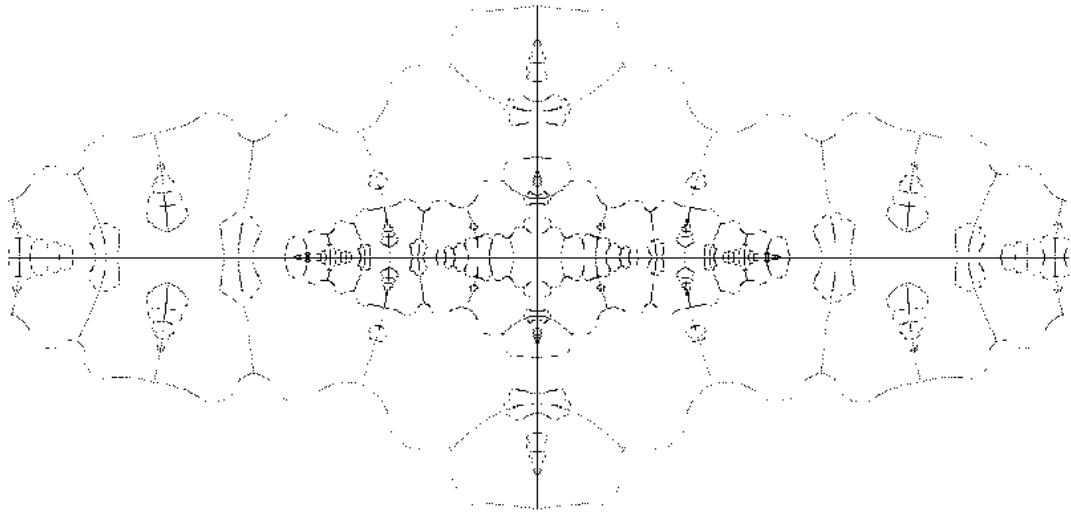}}
\caption{Some pieces of the puzzles ${\cal P}_i$, $i=1,2,3$.} 
\label{puzzle}
\end{figure}

\begin{DEF}
We can define:
\begin{itemize}
\item[$\bullet$]{the puzzle pieces of the puzzle ${\cal P}_0$ of depth $0$ 
are the chess squares $P$, and}
\item[$\bullet$]{the puzzle pieces of the puzzle ${\cal P}_n$ of depth $n$ 
are the pieces $\lambda^n P$.}
\end{itemize}
\end{DEF}
On figure \ref{puzzle}, we have drawn some pieces of the puzzles of depth 
$1$, $2$, and $3$. A lot of pieces are missing (but drawing them with 
a computer is a far too difficult and long work), that is why the puzzles seem
to be disconnected, which is not the case:
remember that the chess squares cover a dense open subset $\widehat  
W\subset\Cm$ which is connected and simply connected.

Lemma \ref{3} enables us to conclude that for every $n\geq 0$, the puzzle 
${\cal P}_n$ covers the Julia set $K(f)$.
On the other hand, it is possible to obtain the following result
which is crucial to show local connectivity.
\begin{LEMME}\label{4}
If $P$ is a chess square, then for all $n\geq 0$, 
the intersection $\lambda^n\overline P\cap K(f)$ is connected.
\end{LEMME}
\dem
Given a chess square $P$, there exists an integer
$k\geq 0$ such that $f_\lambda^k(P)=W_\pm$. We will show lemma \ref{4} by 
induction on the age $k$ of the chess square $P$.

For $k=0$, $\lambda^n \overline W_\pm\cap K(f)$ is connected, because
$f^{2^n}(\lambda^n \overline W_\pm)$ is a half plane, which cuts
$K(f)$ in two parts 
along the real axis, and because $K(f)$ is totally invariant by 
the polynomial-like map $f~:~W\to \Cm_\lambda$.

If $P$ is a chess square of age $k+1\geq 1$, we first notice that
$K(f)\subset W$. Hence, 
\begin{itemize}
\item[$\bullet$]{either $P$ is a copy of $W_\pm$ by rotation of angle 
$2\pi/\l$, and we can easily conclude,}
\item[$\bullet$]{either $P\cap K(f)=\emptyset$.}
\end{itemize}
On the other hand, $f_\lambda(P)=Q$ 
is a chess square of age $k$. Hence, every $\lambda^n \overline Q\cap
K(f)$, $n\in\Nm$, 
is connected.
But, 
$$f(\lambda P)=\left(-\frac{1}{\lambda}\right)^{n-1} f^{2^{n-1}}
(\lambda^{n-1}\lambda P).$$
Hence,
$$f^{2^{n-1}}(\lambda^n P)=(-\lambda)^{n-1}f_\lambda(P)=(-\lambda)^{n-1}Q.$$
We can conclude that $\lambda^n \overline P\cap K(f)$ is connected.\finlem

\begin{LEMME}\label{5}
Let $P$ be a chess square. Then $f(P)=\Hm_\pm$ and  
$f~:~P\to \Hm_\pm$ is an isomorphism.
\end{LEMME}
\dem
The commutative diagram shows that $f$ is the linearizer of $f_\lambda$.
Hence we can write:
$$f(z)=\left(-\frac{1}{\lambda}\right)^n f\circ f_\lambda^n(z).$$
But, if $P$ is a chess square, there exists an integer $n$ such that 
$f_\lambda^n~:~P\to W_\pm$ is an isomorphism. Hence,
$$f(P)=\left(-\frac{1}{\lambda}\right)^n f(W_\pm)=\Hm_\mp.$$
Finally, as $f~:~W_\pm\to \Hm_\mp$ is an isomorphism, we can conclude.
\finlem

We can now show that the puzzle ${\cal P}_n$ 
contains the puzzle ${\cal P}_{n+1}$. In other words, 
we will show the following lemma.
\begin{LEMME}\label{6}
Every puzzle piece $\lambda P\in{\cal P}_1$ is contained in a chess square
$Q\in{\cal P}_0$
\end{LEMME}
This is obvious on figure \ref{puzzle}.
\dem
Let  $z$ be a point in $\lambda P$ and let $Q$ be the chess square which 
contains $z$. Using lemma 5, we can say that $f(Q)=\Hm_\pm$.
As $f(\lambda P)=f_\lambda(P)$ is a chess square, we know that
$$f(\lambda P)\subset f(Q),$$
and as $f~:~Q\to \Hm_\pm$ is an isomorphism, we can deduce that
$$\lambda P\subset Q.$$
\finlem

\begin{figure}[htb]
\centerline{
\psboxscaled{500}{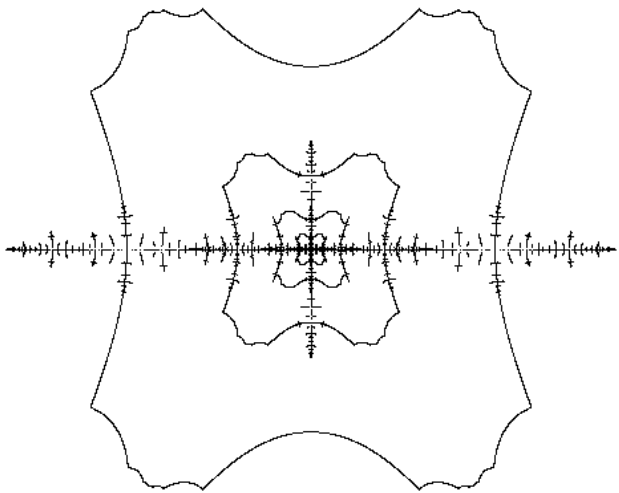}}
\caption{Local connectivity of the Feigenbaum Julia set around the
critical point.} 
\label{local}
\end{figure}

We would like to reprove the following theorem by Jiang and Hu.
\begin{THEO}
The Julia set $K(f)$ is locally connected.
\end{THEO}
Up to now, we can only prove local connectivity around the critical
point $0$ (cf figure \ref{local}):
\begin{itemize}
\item[$\bullet$]{$0\in\lambda^n W\subset \lambda^{n-1} W$, and}
\item[$\bullet$]{$\lambda^n W\cap K(f)$ is connected.}
\end{itemize}

This local connectivity spreads around using classical
arguments of recurrence, but those arguments do not exhibit a nice
basis of connected neighborhoods. We would like to prove that if we
consider the nested sequence of puzzle pieces $P_n(z)\in {\cal P}_n$
which contain a point $z\in K(f)$, then the diameter of the piece
$P_n(z)$ decreases geometrically with $n$.
We just need to show the following conjecture.
\begin{CONJ}
There exists a constant $\varepsilon>0$ such that if $P\in {\cal P}_0$
and $Q\in{\cal P}_0$ are chess squares and $\lambda P\subset Q$, 
then
$$\frac{{\rm diam}(Q)}{{\rm diam}(\lambda P)}\geq 1+\varepsilon.$$
\end{CONJ}
We think it is possible to show this conjecture using some classical
arguments of moduli of annuli. 
Moreover,we think that there is a stronger result, but we have no
idea how to prove it.
\begin{CONJ}
There exists a constant $K\in \Rm$ such that the diameter of every
chess square 
$P\in{\cal P}_0$ is bounded by $K$:
$$\mbox{diam}(P)\leq K.$$
\end{CONJ}

\newcounter{nom}{\setcounter{nom}{1}}

\end{document}